\documentstyle[12pt,epsf]{amsart}

\numberwithin{equation}{subsection}
\newtheorem{theorem}[equation]{Theorem}
\newtheorem{prop}[equation]{Proposition}

\newtheorem{corollary}[equation]{Corollary}
\newtheorem{defn}[equation]{Definition}

\newtheorem{conj}[equation]{Conjecture}

\newcommand{\bd}{\partial}

\newcommand{\reals}{{\Bbb R}}

\newcommand{\TT}{{\Bbb T}}
\newcommand{\KK}{{\Bbb K}}
\newcommand{\SS}{{\Bbb S}}

\newcommand{\Lap}{{\bf \Delta}}

\newcommand{\ph}{\varphi}
\newcommand{\mult}{{\rm mult}}
\newcommand{\cn}{{\rm cn}}

\newcommand{\area}{{\rm Area}}

\begin{document}
\title[Extremal metric for $\lambda_1$ on a Klein bottle]{Extremal
metric for the first eigenvalue on a Klein bottle}
\author[D. Jakobson]{Dmitry Jakobson}
\address{Department of Mathematics and
Statistics, McGill University, 805 Sherbrooke Str. West, Montr\'eal
QC H3A 2K6, Ca\-na\-da.}
\email{jakobson@@math.mcgill.ca}
\author[N. Nadirashvili]{Nikolai Nadirashvili}
\address{Department of Mathematics,
University of Chicago, 5734 University Avenue, Chicago, IL 60637,
USA.}
\email{nicholas@@math.uchicago.edu}
\author[I. Polterovich]{Iosif Polterovich}
\address{D\'e\-par\-te\-ment de math\'ematiques et
de sta\-tistique, Univer\-sit\'e de Mont\-r\'eal
CP 6128 succ Centre-Ville, Mont\-r\'eal QC  H3C 3J7
Canada.}
\email{iossif@@dms.umontreal.ca}


\thanks{The first author was supported by NSERC, FQRNT, Alfred P.
Sloan Foundation fellowship and Dawson fellowship.  The second
author was supported by NSF grant DMS-9971932. The third author
was supported by NSERC and FQRNT}

\begin{abstract}
The first eigenvalue of the Laplacian on a surface can be viewed
as a functional on the space of Riemannian metrics of a given
area. Critical points of this functional are called extremal
metrics. The only known  extremal metrics are a round sphere, a
standard projective plane, a Clifford torus and an equilateral
torus. We construct an extremal metric on a Klein bottle. It is a
metric of revolution, admitting a minimal isometric embedding into
a sphere ${\mathbb S}^4$ by the first  eigenfunctions. Also, this
Klein bottle is a  bipolar surface for the Lawson's
$\tau_{3,1}$-torus. We conjecture that an extremal metric  for the
first eigenvalue on a Klein bottle is unique, and hence it
provides a sharp upper bound for $\lambda_1$ on a Klein bottle of
a given area. We present numerical evidence and prove the first
results towards this conjecture.
\end{abstract}
\maketitle

\section{Introduction and main results}\label{sec:main}
\subsection{Extremal metrics for the first eigenvalue}
Let $M$ be a closed surface of genus $\gamma$ and let $g$  be the
Riemannian metric on $M$. Denote by $\Delta$ the Laplace-Beltrami
operator on $M$, and by $\lambda_1$ the smallest positive
eigenvalue (the fundamental tone) of the Laplacian. How large can
$\lambda_1$ be on such a surface? It was proved in \cite{H},
\cite{YY}, \cite{LY} that
$$\lambda_1\area(M) \le {\rm const}(\gamma),$$
where the constant grows linearly with $\gamma$. However,
for $\gamma \ge 1$ bounds obtained in this way have no reason to
be {\it sharp.} In order to
study sharp upper bounds we recall the following
\begin{defn} A metric $g$ on a surface is called $\lambda_1-maximal$ if
for any metric $\tilde g$ of the same area
$\lambda_1(\tilde g)\le \lambda_1(g)$.
\end{defn}
In other words, a $\lambda_1$-maximal metric is a global maximum
of the functional $\lambda_1 : g \to~\reals$. Consider critical points of this functional.
\begin{defn} {\rm (see \cite{ESI2}).}
An {\it extremal} metric for the first eigenvalue is a critical
point $g_0$ of the functional $\lambda_1: g \to\reals$, i.e. for
any analytic deformation $g_t$ of the Riemannian metric $g_0$ in the class
of metrics of fixed area
$\lambda_1(g_t)\le \lambda_1(g_0) + o(t)$ as $t \to 0$.
\end{defn}
Note that the functional $\lambda_1$ does not have local
minima \cite{ESI2}.

\subsection{Extremal metrics and minimal immersions.}
\label{ex} Only four examples of extremal metrics for the first
eigenvalue are known: (i) standard metric on ${\mathbb S}^2$, (ii)
standard metric on  ${\mathbb RP}^2$, (iii) flat equilateral torus
and (iv) Clifford torus. Moreover, it was proved that there are no
other extremal metrics on  these three surfaces (\cite{MR},
\cite{ESI3}, \cite{ESI2}). Metrics (i)-(iii) are
$\lambda_1$-maximal (\cite{H}, \cite{LY}, \cite{N1}), (iv) is just
a local extremum.

The following remarkable property holds for extremal metrics for
the first eigenvalue. Any surface with an extremal metric admits a
minimal isometric immersion by the first eigenfunctions into  a
round sphere of a certain dimension. In all examples (i)--(iv) the
dimension is equal to $\mult(\lambda_1)-1$, where
$\mult(\lambda_1)$ is the multiplicity of the first eigenvalue.
There is a vast literature on relations between extremal metrics
and minimal immersions, see \cite{Berger}, \cite{LY}, \cite{MR},
\cite{N1}, \cite{ESI3}, \cite{ESI2}.

\subsection{Extremal metric on a Klein bottle}\label{ext}
It is proved in \cite{N1}, that on a Klein bottle there exists a
$\lambda_1$-maximal (and hence an extremal) metric,  which is a
metric  of revolution with $\mult(\lambda_1)=5$.
However, no example of an extremal metric on a Klein bottle has
been known. Our main result is an explicit construction of such a metric.
\begin{theorem}
\label{Klein} A metric of revolution
\begin{equation}
\label{new1}
g_0=\frac{9+(1+8\cos^2 v)^2}{1+8\cos^2 v}\left(du^2
+\frac{dv^2}{1+8\cos^2 v}\right),
\end{equation}
$0\le u< \pi/2, 0\le v < \pi$, is an extremal metric for the first
eigenvalue on a Klein bottle $\KK$. The surface $(\KK, g_0)$
admits a minimal isometric embedding into a sphere ${\mathbb S}^4$
by the first eigenfunctions. The first eigenvalue of the Laplacian
for this metric has multiplicity $5$ and satisfies the equality
\begin{equation}
\label{otvet}
\lambda_1 {\area}(\KK, g_0) = 12\pi E(2\sqrt{2}/3),
\end{equation}
 where
$E(\cdot)$ is a complete elliptic integral of the second kind.
\end{theorem}
\smallskip

\noindent {\bf Remark.} An extremal metric on a Klein bottle must
be a metric  of revolution since any conformal diffeomorphism of an
extremal metric is an isometry (\cite{MR}, \cite{ESI3}), and any
metric on a Klein bottle is conformally equivalent to a flat
metric which is invariant under a natural $\SS^1$-action (see
section 2.1). The condition $\mult(\lambda_1)=5$ follows from the
following argument. It is shown in \cite{ESI3} that
$\mult(\lambda_1)>3$ for an extremal metric on any surface but a
sphere. On the other hand, on a Klein bottle $\mult(\lambda_1)\le
5$ (\cite{N3}), and the case $\mult(\lambda_1)=4$ has been
excluded in \cite{N1}.

\smallskip

We prove Theorem  \ref{Klein} in section 3.

\smallskip

\noindent{\bf Remark.} It is shown in \cite{ESI4} that the
extremal metrics  for the first eigenvalue (i)-(iv)
are also the critical points of the functional $\operatorname{Tr}
e^{-t\Delta}$ (the trace of the heat kernel) at {\em any} time
$t>0$.  Theorem \ref{Klein} shows that this is not always the case:
there
are no critical points of $\operatorname{Tr}
e^{-t\Delta}$ for all $t>0$ on a Klein bottle (\cite{ESI4}).

\subsection{Interpretation in the language of minimal surfaces}
The  Klein bottle $(\KK, g_0)$ constructed in Theorem \ref{Klein}
has the following surprising interpretation in terms of
$\SS^1$-equivariant minimal surfaces in ${\mathbb S}^4$.
Equivariant minimal immersions  into spheres is a classical
subject in minimal surfaces (see \cite{HL}, \cite{Uhlenbeck}). In
particular, $\SS^1$-equivariant minimal immersions of tori and
Klein bottles into ${\mathbb S}^4$ have been studied in \cite{FP}.
\begin{theorem}
\label{bipolar} The surface $(\KK, g_0)$ is  a bipolar surface of
Lawson's $\tau_{3,1}$-torus.
\end{theorem}
In section 4 we prove Theorem \ref{bipolar} and recall the
definitions of Lawson's tori and bipolar surfaces  (see also
\cite{Lawson}).
Interestingly enough, the interpretation of $g_0$ as
a metric on a bipolar surface allows us to
simplify the explicit formula for $g_0$ (cf. (\ref{new1})
and (\ref{metric})).

\subsection{Towards a sharp upper bound for the first eigenvalue}
Combining Theorem \ref{Klein} with the existence of a
$\lambda_1$-maximal metric on a Klein bottle proved in \cite{N1},
we make the following
\begin{conj}\label{conj}
The metric $g_0$ is a unique extremal metric on a Klein bottle, and
in particular it is the $\lambda_1$-maximal metric.
This implies the  following sharp upper bound for the first
eigenvalue on a Klein bottle:
\begin{equation}
\label{con}
\lambda_1(g) {\rm Area}(\KK, g) \le 12\pi E(2\sqrt{2}/3)\approx 13.365 \pi,
\end{equation}
with an equality attained only for $g=g_0$.
\end{conj}
Recall that the estimate of \cite{LY} gives just
$$\lambda_1(g) {\rm Area}(\KK, g) \le 48 \pi.$$
\noindent {\bf Remark.} It is claimed in Theorem 3 in \cite{N1} that
$$\lambda_1(g) {\rm Area}(\KK, g)  \le 8\pi^2/\sqrt{3},$$
the right-hand side being the supremum for $\lambda_1 \area$ on a
torus. However, the proof of this claim is incorrect: it relies on
the  assumption that the first eigenvalue on a Klein bottle is
also the first eigenvalue on the covering torus. Though it is an
eigenvalue on a torus, it might be not the {\it first} eigenvalue.
In particular, for $(\KK, g_0)$ the first eigenvalue is the {\it third}
eigenvalue on the corresponding torus (see Proposition
\ref{Sturm}). Note, however, that indeed $12\pi
E(2\sqrt{2}/3)<8\pi^2/\sqrt{3}.$

\smallskip

In order to prove Conjecture \ref{conj} one has to study the
nonlinear systems of ODEs (\ref{syst01}) or (\ref{syst12}) that
are crucial in the proof of Theorem \ref{Klein}. We need to show
that there are no initial conditions $0<p<1$ except for
$p=\sqrt{3/8}$ (which corresponds to the metric $g_0$) admitting
periodic solutions with the required number of zeros (see
Condition A in section \ref{sec:min}). We discuss numerical
evidence and prove the first results towards Conjecture \ref{conj}
in section 5. However, there are serious difficulties in finding a
rigorous proof of this conjecture, see section 5.7.

\section{A system of ODEs for the extremal metric}

\subsection{Preliminaries}\label{sec:prelim}
We realize the Klein bottle $\KK$ as a fundamental domain in
$\reals^2$ for the group of motions generated by
$(x,y)\to(x+\pi,-y),(x,y)\to(x,y+a)$, where $a>0$ is a conformal
parameter.   $\KK$ has a double cover, the torus ${\Bbb T}^2,$
which is the fundamental domain $\reals^2$ for the group of
motions generated by $(x,y)\to (x+2\pi,y),(x,y)\to(x,y+a)$. The
functions on $\KK$ can thus be thought of as functions on $\TT^2$
satisfying the symmetry condition
\begin{equation}\label{symm}
f(x,y)=f(x+\pi,-y).
\end{equation}
If we expand the functions on $\TT^2$ into Fourier series in $x$,
we can easily see that the functions in $L^2(\TT^2)$ satisfying
\eqref{symm} can be expanded in the series of functions of the
form
\begin{equation}\label{basis1}
\begin{aligned}
\{\phi(y)\sin(2 kx),\phi(y)\cos(2 kx):\\
\phi(y)=\phi(-y),\phi(y+a)=\phi(y)\},
\end{aligned}
\end{equation}
and of the form
\begin{equation}\label{basis2}
\begin{aligned} \{\psi(y)\sin( x(2k+1)),\psi(y)\cos(
x(2k+1)):\\ \psi(y)=-\psi(-y),\psi(y+a)=\psi(y)\},
\end{aligned}
\end{equation}
where $k\in {\Bbb Z}$.

As mentioned in section \ref{ext},
it follows from \cite{N1} that an  extremal metric for the first
eigenvalue on a Klein bottle is necessarily a metric of revolution
and the multiplicity of $\lambda_1$ for this metric is equal to
$5$.

Hence without loss of generality we may assume that our metric is
invariant under the $\SS^1$ action $(x,y)\to(x+t,y)$,
$0\leq t\leq \pi$, and is given by
$\hat g_0=f(y)(dx^2+dy^2)$, where $f(y)=f(y+a)=f(-y)>0$ is
the conformal factor.  The area of the Klein bottle is equal to
$$\area(\KK)=\pi\int_0^a\; f(y)dy.$$
The Laplacian on $\KK$ is given by
\begin{equation}\label{Laplacian}
\Lap=-\frac{1}{f(y)}\left(\frac{\bd^2}{\bd x^2}+\frac{\bd^2}{\bd
y^2}\right)
\end{equation}

Let $\lambda_1$ denote the first nonzero eigenvalue of $\Lap$.
We want to determine the conformal class (i.e. the value of $a$)
that maximizes the product $\lambda_1 \area$.  Since our metric
is rotationally invariant, the operator $\bd/\bd x$ commutes with
$\Lap$, so we can find a joint basis of eigenfunctions of the form
\eqref{basis1} and \eqref{basis2}.
\subsection{First eigenfunctions}
\label{cour} By Courant's nodal domain theorem, any eigenfunction
in the first eigenspace should have exactly two nodal domains. Our
eigenfunctions have the form
\begin{equation}\label{vanish1}
\ph_k(y)\cos(k x),\ \ \ \ph_k(y)\sin(k x),
\end{equation}
where $\ph_k(-y)=(-1)^k\ph_k(y),\ph_k(y+a)=\ph_k(y).$  For $k$
odd, $\ph_k(0)=0$ so $\ph_k$ vanishes at least once.  Also,
$\ph_0$ must vanish at least once since the corresponding
eigenfunction can't have constant sign. Let $\ph_k$ vanish $m_k$
times in the period $[0,a)$.

We can choose the fundamental domain for the Klein bottle to be
the set $X=[y,y+\pi]\times [-a/2,a/2]$, with $y/\pi$ irrational
(to avoid vanishing on the vertical sides), and with the
appropriate boundary identifications.  The nodal set of an
eigenfunction \eqref{vanish1} consists of a grid with $k$ distinct
vertical lines and $m_k$ distinct horizontal lines.  It is easy to
show that such an eigenfunction has at least $k$ nodal domains:
indeed, $k$ vertical lines divide the set $X$ into $k+1$ vertical
strips, and of those only the two boundary strips are glued into
one by side identifications.  Therefore, by Courant's nodal domain
theorem we must have $k\leq 2$.

Substituting into \eqref{Laplacian} and taking into account that
$\mult(\lambda_1)=5$, we conclude that the eigenspace
corresponding to $\lambda$ has a basis of eigenfunctions of the
form
\begin{equation}\label{func:1}
\left\{
\begin{aligned}
\label{eig}
\varphi_0(y), & \qquad \varphi_0(-y)=\varphi_0(y),
\varphi_0''= -\lambda f\varphi_0;\\
\cos(x)\varphi_1(y),&\qquad\varphi_1(-y)=-\varphi_1(y),
\varphi_1''=(1-\lambda f)\varphi_1;\\
\sin(x)\varphi_1(y),&\;\\
\cos(2 x)\varphi_2(y),&\qquad\varphi_2(-y)=\varphi_2(y),
\varphi_2''=(4-\lambda f)\varphi_2;\\
\sin(2 x)\varphi_2(y).&\;
\end{aligned}
\right.
\end{equation}
Here all functions of $y$ are periodic with period $a$ and
$\lambda f$ is an unknown positive function. Since an extremal
metric necessarily admits a minimal isometric immersion into a
sphere (in our case of dimension $4$), we get two more conditions
on the functions $\ph_0$, $\ph_1$ and $\ph_2$ (cf. \cite{N1}):

\begin{equation}\label{crit1}
\varphi_0^2+\varphi_1^2+\varphi_2^2=1.
\end{equation}

\begin{equation}\label{crit2}
(\varphi_0')^2+(\varphi_1')^2+(\varphi_2')^2=\varphi_1^2+4\varphi_2^2
=\lambda f/2.
\end{equation}

We can now substitute for $\lambda f$ in the second and the
third equations in \eqref{func:1}, getting the following system of
second order equations for $\varphi_1$ and $\varphi_2$ (where
$\lambda f$ has been eliminated):
\begin{equation}\label{func:2}
\left\{
\begin{aligned}
\varphi_1''=(1-2(\varphi_1^2+4\varphi_2^2))\varphi_1;\\
\varphi_2''=(4-2(\varphi_1^2+4\varphi_2^2))\varphi_2.
\end{aligned}
\right.
\end{equation}

\subsection{Zeros of the first eigenfunctions}
\label{sec:min} We use the Courant's nodal domain theorem once
again (see previous section) to get a condition on the number of
zeros of $\ph_0, \ph_1$ and $\ph_2$. Since each of these functions
is a non-trivial solution of a second order differential equation
(\ref{func:1}), it is impossible that $\ph_k$ and $\ph_k'$ vanish
simultaneously for $k=0,1,2$. Periodicity then implies that the
number of zeros $m_k$ for any $\ph_k$ is an even number. Recalling
that each eigenfunction has exactly two nodal domains, and taking
into account boundary identifications as in the previous section,
we get $m_0=m_1=2$ and $m_2=0$.

{\bf Condition A (zeros).} $\ph_0$ and $\ph_1$ should have exactly
{\em two} zeros in the period, while $\ph_2$ {\em should not
vanish}.

\subsection{First integrals}\label{sec:int}
It is straightforward to check that the following expressions are
the first integrals for the system \eqref{func:1}, \eqref{crit1}
(cf. \cite{Uhlenbeck}):
\begin{equation}\label{klein:integr2}
\left\{
\begin{aligned}
E_0:=\varphi_0^2+(\varphi_0\varphi_1'-
\varphi_1\varphi_0')^2+(\varphi_0\varphi_2'-\varphi_2\varphi_0')^2/4,\\
E_1:=\varphi_1^2+(\varphi_1\varphi_2'-
\varphi_2\varphi_1')^2/3-(\varphi_1\varphi_0'-\varphi_0\varphi_1')^2,\\
E_2:=\varphi_2^2-(\varphi_2\varphi_0'-
\varphi_0\varphi_2')^2/4-(\varphi_2\varphi_1'-\varphi_1\varphi_2')^2/3.
\end{aligned}
\right.
\end{equation}
In the verification of this fact, one uses \eqref{crit1} and its
consequence $\ph_0\ph_0'+\ph_1\ph_1'+\ph_2\ph_2'=0$. In fact, all
these integrals are equivalent: one can show that
$E_0+E_1+E_2=1=E_0+3E_1/4$, $E_2=-E_1/4$.  Hence $E_j$'s define
just one independent first integral. We make use of the different
expressions (\ref{klein:integr2}) in section 5.

Let us evaluate $E_1$ at $y=0$.
\begin{equation}
\label{init1}
\varphi_1(0)=0=\varphi_0'(0)=\varphi_2'(0)
\end{equation}
It follows that $\varphi_0(0)^2+\varphi_2(0)^2=1$ and that
\begin{equation}
\label{init2}
\varphi_1'(0)^2=4\varphi_2(0)^2.
\end{equation}
Substituting into
the expression for $E_1$ we find that
\begin{equation}\label{int:incond}
E_1=\frac{4}{3}\varphi_2(0)^2(4\varphi_2(0)^2-3).
\end{equation}
\subsection*{Remark}
Alternatively, one can start with the system \eqref{func:2} (or
similar systems involving just $\ph_0,\ph_1$, or just
$\ph_0,\ph_2$) and deduce the following equivalent expressions for
the first integrals \eqref{klein:integr2}:
\begin{equation}\label{klein:integr1}
\left\{
\begin{aligned}
(\varphi_1')^2+4(\varphi_2')^2+(\varphi_1^2+
4\varphi_2^2)^2-\varphi_1^2-16\varphi_2^2:=\kappa_0,\\
(\varphi_0')^2-3(\varphi_2')^2+2\varphi_0^2+
6\varphi_2^2-(\varphi_0^2-3\varphi_2^2)^2:=\kappa_1,\\
4(\varphi_0')^2+3(\varphi_1')^2+32\varphi_0^2+
21\varphi_1^2-(4\varphi_0^2+3\varphi_1^2)^2:=\kappa_2.
\end{aligned}
\right.
\end{equation}
One can show that $\kappa_2-3\kappa_0-4\kappa_1=12$ and that
$\kappa_0+\kappa_1=1$, so $\kappa_0+\kappa_2=16$. One can also
show that $E_1=\frac{1}{3}\kappa_0$. For certain applications it
is more convenient to use the  expressions \eqref{klein:integr1}
rather than \eqref{klein:integr2}; however, we shall not use the
expressions \eqref{klein:integr1} in this paper.

\section{Proof of Theorem \ref{Klein}}\label{sec:period}
\subsection{A system for $\ph_0$ and $\ph_1$}
The initial conditions in system (\ref{func:2}) can be
parametrized as follows:
\begin{equation}\label{syst12}
\left\{
\begin{aligned}
\varphi_1''=(1-2\ph_1^2-8\ph_2^2)\varphi_1\\
\varphi_2''=(4-2\ph_1^2-8\ph_2^2)\varphi_2\\
\varphi_1(0)=0, \,\, \varphi_1'(0) = 2p\\
\varphi_2(0)=p, \,\, \varphi_2'(0)=0,
\end{aligned}
\right.
\end{equation}
where $0 \le p\le 1$ is a parameter of the system. Moreover,
$\varphi_0$ and $\varphi_1$ both have two zeros on the
period, while $\varphi_2$ has constant sign.

The corresponding system for the functions $\ph_0, \ph_1$ reads:
\begin{equation}\label{syst01}
\left\{
\begin{aligned}
\varphi_0''=(8\varphi_0^2+6\varphi_1^2-8)\varphi_0\\
\varphi_1''=(8\varphi_0^2+6\varphi_1^2-7)\varphi_1\\
\varphi_1(0)=0, \,\, \varphi_1'(0) = 2p\\
\varphi_0(0)=\sqrt{1-p^2}, \,\, \varphi_0'(0)=0,
\end{aligned}
\right.
\end{equation}

Note that in (\ref{syst12}) and (\ref{syst01}) initial conditions
are determined by (\ref{init2}) modulo signs.
However, changing the signs of the initial conditions may
only result in changing the signs of the solutions
(in other words, we will get the same eigenfunctions possibly
multiplied by $-1$). Therefore, we may consider only
non-negative initial conditions in (\ref{syst12}) and
(\ref{syst01}).

\subsection{Solution for $p=\sqrt{3/8}$}
\label{sol}
Our objective is to find values of $p$ such that the system has
periodic solutions satisfying Condition A, namely that both
$\ph_0$ and $\ph_1$ have exactly two zeros on the period. We find
a candidate from a numerical experiment: $p=\sqrt{3/8}$. Note that
this value of $p$ is exactly the minimum of the first integral
$E_1$ and hence as follows from \cite{FP} it corresponds to a
periodic solution. We discuss this in more detail in section 4.

Set $p=\sqrt{3/8}$. Let us look for $\ph_0$ and $\ph_1$ in the
following form:
\begin{equation}\label{zamena}
\ph_0(y)=\sqrt{\frac{5}{8}}\cos \theta (y), \,\,
\ph_1(y)=\frac{1}{\sqrt{2}} \sin \theta(y), \,\, \theta(0)=0, \,\,
\theta'(0)=\sqrt{3}
\end{equation}
Such a change of variables is also motivated by numerical
experiments, suggesting that
\begin{equation}\label{rel}
2 \ph_1^2 +8/5 \ph_0^2=1.
\end{equation}
Initial conditions for $\theta$ are prescribed by the initial conditions
for $\ph_1, \ph_0$.

Of course, in principle, such an ansatz could make our system
overdetermined: note that instead of two variables $\ph_0, \ph_1$
we now have one variable $\theta$. However, as shown below,
for this particular choice of constants this does not happen.

Indeed, we have:
$$8 \ph_0^2 + 6 \ph_1^2 = 5 \cos^2 \theta +
3\sin^2 \theta=5-2\sin^2\theta,$$
and hence (\ref{syst01}) can be rewritten as
\begin{equation}\label{theta1}
\left\{
\begin{aligned}
(\theta')^2 -\theta''
\frac{\cos \theta}{\sin\theta} = 2 + 2 \sin^2 \theta\\
(\theta')^2+\theta''
\frac{\sin \theta}{\cos \theta} = 3 + 2 \sin^2 \theta\\
\theta(0)=0, \,\, \theta'(0) = \sqrt{3}.
\end{aligned}
\right.
\end{equation}
Subtracting the second equation from the first we get
\begin{equation}
\label{theta2}
\theta'' = \sin\theta \cos\theta=\frac{1}{2} \sin 2\theta
\end{equation}

Multiplying by $\theta'$ and integrating gives
\begin{equation}
\label{theta3}
(\theta')^2 = 3 + \sin^2\theta.
\end{equation}
Exactly the same equation one gets if (\ref{theta2}) is substituted
into (\ref{theta1}) and hence the whole system yields to
(\ref{theta3}) with an initial condition $\theta(0)=0$, implying
$$y=\frac{1}{2}\int_0^\theta
\frac{d\theta}{\sqrt{1-\frac{1}{4}\cos^2\theta}}$$ From this
equation we can deduce periodicity conditions. The functions
$\ph_0, \ph_1$ are periodic in $\theta$ with the period $2\pi$.
Hence, the period $a$ is equal to
$$
\begin{aligned}
\frac{1}{2}\int_0^{2\pi}
\frac{d\theta}{\sqrt{1-\frac{1}{4}\cos^2\theta}}=
2\int_0^{\pi/2}\frac{d\theta}{\sqrt{1-\frac{1}{4}\cos^2\theta}}=\\
2\int_0^{\pi/2}\frac{d\theta}{\sqrt{1-\frac{1}{4}\sin^2\theta}}=
2K(1/2),
\end{aligned}
$$
where $K(\cdot)$ is a complete elliptic integral of the
first kind. Hence, $a=2 K(1/2)$.

Let us now compute $\lambda \area({\KK})$ for this metric (even
without computing the metric explicitly -- this will be done in
the next section). Taking into account (\ref{crit2}) we have:
$$
\lambda {\rm Area}(\KK)=\lambda \pi \int_0^a f(y) dy = 2\pi
\int_0^{2K(1/2)} 4- 3\ph_1^2(y)-4\ph_0^2(y) dy=$$
$$
2\pi \int_0^{2K(1/2)}4 - 5/2 \cos^2 \theta - 3/2 \sin^2\theta dy=
$$
$$
2\pi \int_0^{2\pi} (5/2 - \cos^2 \theta) y'(\theta) d\theta=
2\pi \int_0^{2\pi} \frac{5/2 - \cos^2\theta}{\sqrt{4-\cos^2 \theta}}
d\theta =
$$
$$
2\pi \left(\int_0^{2\pi} \sqrt{4 - \cos^2\theta}d\theta -
3/2 \int_0^{2\pi} \frac{d\theta}{\sqrt{4-\cos^2\theta}}\right)=
$$
$$
2\pi\left(8\int_0^{\pi/2} \sqrt{1-1/4 \cos^2 \theta}d\theta -
3 \int_0^{\pi/2} \frac{d\theta}{\sqrt{1 - 1/4 \cos^2\theta}}\right)=
$$
$$
2\pi (8 E(1/2) - 3 K(1/2)) = 12\pi E(2\sqrt{2}/3).
$$
The last equality follows from an identity relating the complete
elliptic integrals of the first and the second kind (see
\cite{Erd}, p. 319). This proves the assertion (\ref{otvet}) in
Theorem \ref{Klein} up to the fact that $\lambda$ is the first
eigenvalue (see section \ref{firstei}).

\subsection{The eigenfunctions}
In this section we find explicitly the eigenfunctions
corresponding to the value $p=\sqrt{3/8}$, and the corresponding
metric $\hat g_0$. We do this using the relation (\ref{rel})
between $\ph_0$ and $\ph_1$, which implies a similar equation for
$\ph_2$ and $\ph_1$: $4 \ph_2^2 - \ph_1^2 = 3/2$. This allows to
transform our system into three separate equations on $\ph_0,
\ph_1, \ph_2$:
$$\ph_0'' = 16/5 \ph_0^3 - 5\ph_0,\,\,
\ph_1''=-2\ph_1 - 4\ph_1^3,\,\,
\ph_2'' = 7 \ph_2 - 16\ph_2^3.$$
We then reduce them to first order equation:
$$(\ph_0')^2 = 8/5 \ph_0^4 - 5\ph_0^2 + 20/8,$$
$$(\ph_1')^2 = -2\ph_1^2 - 2\ph_1^4 + 3/2$$
$$(\ph_2')^2 = 7 \ph_2^2 -8\ph_2^4 - 3/2.$$
Each of these equation can be solved in terms of
elliptic functions (see also \cite{WW}, section 20.6).
Finally we get:
\begin{equation}
\label{phi0}
\ph_0(y)=\sqrt{5/8}\left(1-
\frac{3}{2 \wp(y;73/12, -595/216)-\frac{1}{6}}\right),
\end{equation}
$$\ph_1(y)=\frac{1}{\sqrt{2}}\left(-1 + \frac{2}
{\wp(y+\frac{K(1/2)}{2}; -8/3, 28/27) +\frac{2}{3}}\right),
$$
$$
\ph_2(y)=\sqrt{3/8} +\frac{1}{4}\left(\frac{\sqrt{3/2}}
{\wp(y;193/12, 2681/216) + \frac{11}{12}}\right),
$$
where $\wp(y;\gamma_1,\gamma_2)$ is a Weierstrass $\wp$-function with
invariants  $\gamma_1,\gamma_2$.

It can be checked directly (analytically or using {\it Mathematica})
that $\ph_0, \ph_1, \ph_2$  satisfy the Condition A of
section 2.2.

Using (\ref{rel}) we find that the normalized
metric $\hat g_0=\lambda f(y)(dx^2+dy^2)$
(though it differs by a normalization factor $\lambda$
from the metric defined in the beginning of section 2.2,
we denote it also $\hat g_0$) is given by
\begin{equation}
\label{metric}
\lambda f(y) = 2 (\ph_1^2(y) + 4 \ph_2^2(y))=
5 - \frac{16}{5}\ph_0^2(y),
\end{equation}
The metric $\hat g_0$ is conformally equivalent to a flat metric
on $\KK$ corresponding to the lattice $(x,y)\to (x+\pi, -y)$,
$(x,y)\to (x,y+2K(1/2))$. It still remains to show that this
metric coincides (up to a dilatation) with the metric $g_0$
defined by (\ref{new1}). We postpone this until section~4.
\subsection{Why $\lambda$ is the first eigenvalue?}
\label{firstei}
To complete the proof of Theorem \ref{Klein} we need to show that the eigenvalue $\lambda=1$
of the normalized metric $\hat g_0$ (\ref{metric}) is the {\it first} eigenvalue
of the Laplacian (the eigenvalue equals $1$ due to the choice
of normalization (\ref{metric})).
Note that though Condition A is a necessary condition for the first eigenfunctions,
apriori it is not sufficient.
\begin{prop}
\label{Sturm} The eigenvalue $\lambda=1$ corresponding to the
eigenfunctions
$\{\ph_0(y)$, $\ph_1(y) \cos x$, $\ph_1(y) \sin x$,  $\ph_2(y) \cos 2x$, $\ph_2(y) \sin 2x)\}$,
is the first eigenvalue of the Laplacian on a Klein bottle $(\KK, \hat g_0)$.
\end{prop}
\noindent{\bf Proof.} We prove this proposition with the help
of oscillation theorems of Haupt and Sturm (see \cite{Cod}, \cite{Bee}).
As was mentioned in section 2.2, due to Courant's nodal domain theorem
the first eigenvalue on a Klein bottle of revolution
can be obtained only from one of the three periodic Sturm-Liouville equations
(\ref{func:1}). We need to show that none of the these equations has an eigenfunction
corresponding to an eigenvalue $\lambda < 1$ and satisfying the Condition A as well as
the parity conditions (we are interested only in even eigenfunctions of the
first and the third equation, and only in odd eigenfunctions of the second equation).

Equations  (\ref{func:1}) are subject to
a theorem of Haupt, stating that each eigenvalue problem has a
sequence of eigenfunctions
$$\psi_0, \psi_1, \psi_2, \dots, \psi_{2n-1}, \psi_{2n}, \cdots\,
$$
such that $\psi_0$ does not have zeros and $\psi_{2n-1}, \psi_{2n}$
have exactly $2n$ zeros.
Let us study the equations (\ref{eig})
for $\ph_0, \ph_1$ and $\ph_2$ separately.

The easiest case is $\ph_2$ -- it has no zeros, and hence
corresponds to the smallest eigenvalue of the Sturm-Liouville
problem.

To handle $\ph_1$, we note that since it is odd, it should have
zeros and therefore $0$-th eigenvalue
for this problem is automatically out of question (indeed, this
eigenvalue
$\tilde \lambda\approx 0.2517$ is the first eigenvalue
on the corresponding torus  which covers our Klein bottle).
Since $\ph_1$
has exactly two zeros, it is either $\psi_1$ or $\psi_2$.
Assume there exists another odd solution $\tilde \ph_1$ of the
eigenvalue problem with exactly $2$ zeros, and the corresponding
$\tilde \lambda < \lambda$. Then by Sturm theorem, between each
zero of $\tilde \ph_1$ there should be zeros of $\ph_1$, but since
both of them vanish at $0$, this will mean that $\ph_1$ should
have at least $3$ zeros, while it has only $2$, and we get a
contradiction.
Indeed, numerically one can see that there exists another (but even)
eigenfunction with $2$ zeros. For the normalized problem it
corresponds to $\tilde \lambda \approx 1.31$.

A similar argument works for $\ph_0$. Note that in this case the
$0$-th eigenfunction is just a constant, and hence is also not
relevant. Similarly, $\ph_0$ is either the first or the second
eigenvalue. Assume there is another even eigenfunction with
exactly $2$ zeros corresponding to some $\tilde \lambda <
\lambda$. Note that since it is even and periodic with period $1$
it should be symmetric with respect to the mid-point of the period
$y=1/2$ (indeed, $\ph_0(x)=\frac{\ph_0(x)+\ph_0(1-x)}{2}$), in
particular its zeros have this symmetry --- as do the zeros of
$\ph_0$. On the other hand, due to Sturm theorem, between each
zero of $\tilde \ph_0$ there should be a zero of $\ph_0$, or in
other words zeros of $\ph_0$ and $\tilde \ph_0$ should interlace,
but this contradicts the fact that they have the symmetry property
(symmetry implies that two zeros of one eigenfunction are
between two zeros of the other). Numerically one can
observe that in reality $\tilde \ph_0$ is an odd function
corresponding to $\tilde \lambda \approx 0.7768$. This completes
the proof of the proposition. \qed

\subsection{End of the proof of Theorem \ref{Klein}}
Let us summarize the results of sections 3.1 -- 3.4. We have
constructed a metric of revolution $\hat g_0$ on a Klein bottle,
admitting an isometric embedding into ${\mathbb S}^4$ by the first
eigenfunctions. The first eigenvalue for this metric satisfies the
equality
$$\lambda_1 {\area}(\KK, \hat g_0) = 12\pi E(2\sqrt{2}/3),$$
where $E(\cdot)$ is a complete elliptic integral of the second
kind. Hence to complete the proof we just need to show that the
metric $\hat g_0$ is indeed an extremal metric for the first
eigenvalue. We use Proposition 1.1 of \cite{ESI2}, implying that
if the isometric immersion is given by a complete set of the first
eigenfunctions (i.e. if the eigenfunctions form a basis of the
corresponding eigenspace), then the metric is extremal for
$\lambda_1$. This is clearly the case for us, since we have used
all five eigenfunctions to construct the immersion, and
$5$ is the maximal possible multiplicity for the first eigenvalue
on a Klein bottle. It follows
from the explicit formulas for the eigenfunctions in section 3.3
that it is in fact an embedding.
To complete the proof of Theorem \ref{Klein} it remains to show
that the metric $\hat g_0$ (\ref{metric}) coincides up to a dilatation
with  the metric
$g_0$ (\ref{new1}). This is done in section \ref{bip}
while proving Theorem \ref{bipolar}. \qed

\section{The extremal metric and $\SS^1$-equivariant immersions}
\subsection{Bipolar surfaces and Lawson tori}
In this section we follow \cite{Lawson} and \cite{Ken}.
Let $\mu:M\to \SS^3\subset {\Bbb R}^4$ be a minimal isometric
immersion  of a surface $M$ into $\SS^3$. A Gauss map
$\mu^*:M\to \SS^3$ is defined pointwise as the image
of the unit normal in $\SS^3$ translated to the origin in
${\Bbb R}^4$. The image $\mu^*(M)$ is called a {\it polar variety}.

Let $\tilde \mu=\mu\wedge\mu^*$~--- the
exterior product of $\mu$ and $\mu^*$. It is a vector in
$\wedge^2 {\Bbb R}^4={\Bbb R}^6$,
and one can check (\cite{Lawson})
that it defines a minimal immersion of $M$ into $\SS^5$:
$\tilde \mu: M \to \SS^5\subset {\Bbb R}^6$.
The minimal surface
$\tilde M = \tilde \mu (M)$ in $\SS^5$ is called a
{\it bipolar surface} to $M$. The metric on $\tilde M$ is given by
$ds^2 = (2-\kappa)d\sigma^2$, where
$d\sigma^2$ is the metric on $M$ and
$\kappa$ is the Gaussian curvature on $M$ (\cite{Lawson}).

Let $M=\tau_{m,k}$ ($m\ge k\ge 1$)
be a Lawson's torus, that is a minimal torus
defined by a doubly periodic immersion $\mu:\reals^2\to \SS^3$,
\begin{equation}
\label{ltorus}
\mu(u,v)=(\cos mu \cos v, \sin mu \cos v, \cos ku \sin v,
\sin ku \sin v).
\end{equation}
One may check that the bipolar surface for $\tau_{m,k}$ is a
minimal torus or a minimal Klein bottle in $\SS^4$
(\cite{Lawson}). The metric on a bipolar surface for $\tau_{m,k}$
is given by (see \cite{Ken}):
\begin{equation*}
\label{bipmet}
ds^2=\frac{(k^2 + (m^2-k^2) \cos^2 v)^2+m^2 k^2}{k^2 + (m^2-k^2)
\cos^2 v}
\left(du^2+\frac{dv^2}{k^2 + (m^2-k^2) \cos^2 v}\right).
\end{equation*}
\subsection{Proof of Theorem \ref{bipolar}}
\label{bip} Note that for $m=3$, $k=1$ the metric above is exactly
the metric $g_0$ (\ref{new1}). Let us check that $(\KK, \hat g_0)$
defines a bipolar surface to the $\tau_{3,1}$-torus. One should
verify that $\hat g_0$ coincides (up to a dilatation) with
(\ref{new1}) (the rest is straightforward). This result can be
deduced from the arguments of (\cite{FP}). Indeed, due to
(\ref{eig}) the metric $g_0$ determines a $(2,1)$-equivariant
minimal immersion in $\SS^4$. Moreover, the first integral $E_1$
($H_2$ in the notations of \cite{FP}) achieves its minimum for
$p=\sqrt{3/8}$, see section \ref{sol}. Hence, as mentioned in
(\cite{FP}, p. 274), this metric defines a bipolar surface for the
Lawson's torus in $\SS^3$ corresponding to $(2+1,2-1)$ circular
action, that is exactly the torus $\tau_{3,1}$. The equation
(\ref{rel}) defining the extremal metric is equivalent to equation
(11) in \cite{FP} by setting $z:=\ph_2$, $w:=\ph_1$. Therefore,
metric $\hat g_0$ indeed coincides with $\hat g_0$ up to a
rescaling. This completes the proof of Theorem \ref{bipolar} and
also finishes the last step of the proof of Theorem \ref{Klein}.
\qed

\smallskip

\noindent{\bf Remark.} In fact, it can be verified directly that
$\hat g_0 = 2 g_0$. It is a lengthy calculation in elliptic
functions. Indeed, set $x=u$ and
$$z(v)=\int_0^v \frac{dv}{\sqrt{1+8\cos^2 v}}=
\frac{1}{3}\int_0^v \frac{dv}{\sqrt{1-\frac{8}{9}\sin^2 v}}$$
Then in  the $(x,z)$ variables metric (\ref{new1})
becomes conformal. Note also that the relation above implies
$\cos v = \cn (3z, 2\sqrt{2}/3)$,
the corresponding Jacobi elliptic function.
Set $y=2z+\frac{K(1/2)}{2}$ (note that $2 K(1/2)=\frac{4}{3}K(2\sqrt{2}/3)$
is the period of $\cn (3z, 2\sqrt{2}/3)$). Taking into account (\ref{phi0})
and (\ref{metric}) we arrive to the following identity that it suffices
to check:
$$
\frac{(1 + 8 \cn^2 (3z, \frac{2\sqrt{2}}{3})^2+9)}{1 + 8 \cn^2
(3z, \frac{2\sqrt{2}}{3})^2}=
10-\left(\frac{24\wp(y;\frac{73}{12},-\frac{595}{216})-38}
{12\wp(y;\frac{73}{12},-\frac{595}{216})-1}\right)^2.
$$
The clue to this identity is the following relation between the
Jacobi and the Weierstrass elliptic functions (see \cite{Erd},
13.16.5):
$$
\cn^2\left(3z,\frac{2\sqrt{2}}{3}\right)=\frac{12
\wp(2z;\frac{73}{12},-\frac{595}{216})-10}
{12 \wp(2z;\frac{73}{12},-\frac{595}{216})+17}.
$$
The remainder of the argument is a rather straightforward
application of formulas from section 13.13 of  \cite{Erd}.

\section{Towards a sharp upper bound for the first eigenvalue}
\subsection{Two intervals of the parameter}
The aim of section 5 is to present numerical evidence for
Conjecture \ref{conj} and to prove the first result in that
direction (Theorem \ref{thm:ruleout1}). Our ultimate goal is to
show that there are no extremal metrics corresponding to the
values of the parameter $0< p < 1$ except for $p=\sqrt{3/8}$. It
turns out that the dynamics of the solutions differs for $0< p <
\sqrt{3}/2$ and $\sqrt{3}/2\leq p < 1$. We study these two
intervals separately. In the latter case we prove the absence of
extremal metrics (sections 5.2 and 5.3). For $0<p<\sqrt{3}/2$ we
present a purely numerical argument (section 5.6) and explain the
nature of difficulties in proving Conjecture \ref{conj} (section
5.7).

Initial conditions of (\ref{syst01}) and (\ref{syst12}) are
parametrized by values of $0<p<1$. We shall be using first
integrals (\ref{klein:integr2}):
$$
E_1=(4/3)p^2(4p^2-3),\ E_2=(-1/3)p^2(4p^2-3),\ E_0=1-p^2(4p^2-3).
$$
The periodic solution corresponds to $p^2=3/8$ which is the
minimum of $E_1$. We want to show that there are no other periodic
solutions satisfying Condition A.

The value $p^2=3/4,E_1=E_2=0,E_0=1$ corresponds to a separatrix of
some sort, the behavior of the solutions changes (section 5.3).

\subsection{Ruling out the interval $1> p > \sqrt{3}/2$}\label{sec:sphere1}
In this section we show that there are no periodic solutions of
\eqref{syst12} satisfying Condition A, if $\sqrt{3}/2<p<1$, and
hence $E_1>0,E_2<0,0<E_0<1$.
\begin{theorem}\label{thm:ruleout1}
Assume that $\sqrt{3}/2<p<1$.  Then the system \eqref{func:1}
doesn't have periodic solutions satisfying Condition A.
\end{theorem}
We shall prove Theorem
\ref{thm:ruleout1} by showing that the solutions of the system
(\ref{syst12}) ``rotate'' around the origin in the
$(\ph_1,\ph_2)$-plane.  In other words, if we introduce polar
coordinates in \eqref{syst12}, the angle will be monotone
increasing.  This fact implies that the function $\varphi_2$
vanishes at {\em some} point on any periodic orbit, contradicting
the condition A in section \ref{sec:min}. We use the condition
\eqref{crit1} to introduce spherical coordinates in the system
\eqref{func:1} and use the integrals $E_1$ and $E_2$ to rule out
the initial conditions $\sqrt{3}/2<p<1$.

We introduce the following spherical coordinates in
\eqref{func:1}:
\begin{equation}\label{sphere:1}
\left\{
\begin{aligned}
\varphi_0 &=\cos\psi,\\
\varphi_1 &=\sin\psi\sin\theta,\\
\varphi_2 &=\sin\psi\cos\theta.
\end{aligned}
\right.
\end{equation}
Taking into account parity conditions in \eqref{func:1}, we find
that $\psi$ is an even function, while $\theta$ is an odd
function.

Differentiating once, we find that
\begin{equation}\label{sphere:der}
\begin{aligned}
\varphi_0'=-\sin\psi\cdot\psi',
\varphi_1'=\cos\psi\cdot\psi'\sin\theta+\sin\psi\cos\theta\cdot\theta',\\
\varphi_2'=\cos\psi\cdot\psi'\cos\theta-\sin\psi\sin\theta\cdot\theta'.
\end{aligned}
\end{equation}
It is easy to see that $\psi$ and $\theta$ satisfy the following
initial conditions:
\begin{equation}\label{sphere:IC}
\begin{aligned}
\theta(0)=0,\theta'(0)=2,\\
\psi'(0)=0,\psi(0)=\arcsin{p}.
\end{aligned}
\end{equation}

We next express the integrals $E_0,E_1,E_2$ in terms of
$\theta,\psi$ and their derivatives.  An elementary calculation
using \eqref{sphere:1} and \eqref{sphere:der} gives the following
identities:
\begin{equation}\label{sphere:7}
\left\{
\begin{aligned}
\ph_0\ph_1'-\ph_1\ph_0' &=\psi'\sin\theta+\frac{\sin(2\psi)}{2}
\cos\theta\cdot\theta';\\
\ph_0\ph_2'-\ph_2\ph_0' &=\psi'\cos\theta-\frac{\sin(2\psi)}{2}
\sin\theta\cdot\theta';\\
\ph_1\ph_2'-\ph_2\ph_1' &=-\sin^2\psi\cdot\theta'.
\end{aligned}
\right.
\end{equation}

We now substitute \eqref{sphere:1} and \eqref{sphere:7} into
\eqref{klein:integr2}. We obtain the following identities  for
$E_1$ and $E_2$:
\begin{equation}\label{sphere:E1}
\begin{aligned}
E_1=\sin^2\theta(\sin^2\psi-(\psi')^2)-
\frac{\psi'\theta'\sin(2\psi)\sin(2\theta)}{2}\\
+\frac{(\theta')^2\sin^4\psi}{3}-
\left(\frac{\theta'\cos\theta\sin(2\psi)}{2}\right)^2.
\end{aligned}
\end{equation}

\begin{equation}\label{sphere:E2}
\begin{aligned}
E_2=\cos^2\theta(\sin^2\psi-(\psi')^2/4)+
\frac{\psi'\theta'\sin(2\psi)\sin(2\theta)}{8}\\
-\frac{(\theta')^2\sin^4\psi}{3}-
\left(\frac{\theta'\sin\theta\sin(2\psi)}{4}\right)^2.
\end{aligned}
\end{equation}

If we now assume that $\theta'=0$, the two expressions simplify to
\begin{equation}\label{E1E2:theta0}
\left\{
\begin{aligned}
E_1=\sin^2\theta(\sin^2\psi-(\psi')^2),\\
E_2=\cos^2\theta(\sin^2\psi-(\psi')^2/4).
\end{aligned}
\right.
\end{equation}
We now prove the main theorem of this section.

\noindent{\bf Proof of Theorem \ref{thm:ruleout1}.} Assume that
there exists a periodic orbit such that $\ph_2\neq 0$ (this is one
of the requirements of
 condition A of section \ref{sec:min}). Then this orbit
has a point satisfying $\theta'=0$, since $\theta$ is the angle in
polar coordinates in $(\ph_1,\ph_2)$-plane, and the orbit is a
compact set lying in the upper half-plane.  Now, if
$\sqrt{3}/2<p<1$, we have $E_1>0,E_2<0$.  We evaluate those
integrals at a point where $\theta'=0$.  Substituting into
\eqref{E1E2:theta0}, we see that
\begin{equation}\label{E1E2:contradiction}
\left\{
\begin{aligned}
\sin^2\psi-(\psi')^2>0,\\
\sin^2\psi-(\psi')^2/4<0.
\end{aligned}
\right.
\end{equation}
It is clear that \eqref{E1E2:contradiction} leads to a
contradiction.  This finishes the proof of Theorem
\ref{thm:ruleout1}. \qed

\subsection{The separatrix $p=\sqrt{3}/2$}
If $p=\sqrt{3}/2$, the first integral $E_1$ vanishes.
The solutions are given explicitly by
the following formulas:
\begin{equation}\label{sqrt32:formula}
\ph_0(y)=(3\cos(\theta(y))-1)/4,
\ph_1(y)=\sqrt{3}\sin(\theta(y))/2,
\end{equation}
where
$$
\theta(y)=\pi - 4 \,{\rm arccot}(e^y).
$$
One can observe that these solutions are not periodic. As $y\to
\infty$, we get the upper half of an ellipse in the plane $\ph_0,
\ph_1$.

In fact, since for $p=\sqrt{3}/2$ the integral $E_1=0$,
if there existed a corresponding minimal isometric immersion
of a Klein bottle into ${\mathbb S}^4$
it would be superminimal (see \cite{FP}). However,
as indicated in the appendix of \cite{MR}, the only
superminimal surface immersed into ${\mathbb S}^4$ by the first
eigenfunctions is the standard sphere (this result is attributed
to N. Ejiri).

\subsection{Dynamics in the $(\ph_0,\ph_1)$-plane for
$0\le p<\sqrt{3}/2$} \label{sec:elem}
We are left to check that the solution for $p=\sqrt{3/8}$ is the only
one in the
interval $0<p<\sqrt{3}/2$.  For the first integrals, this interval
corresponds to
$E_1<0,E_2>0,1<E_0<25/16$.
\begin{prop}\label{prop:1}
For $0<p<\sqrt{3}/2$, functions $\ph_2$ and
$\ph_1\ph_0'-\ph_0\ph_1'$ don't vanish.  Moreover, the function
$\varphi_2$ is bounded away from $\pm 1$.
\end{prop}

{\bf Proof.} We recall from \eqref{klein:integr2} that
\begin{equation}\label{E2:expr}
E_2=\varphi_2^2-(\varphi_2\varphi_0'-
\varphi_0\varphi_2')^2/4-(\varphi_2\varphi_1'-\varphi_1\varphi_2')^2/3.
\end{equation}
Since $E_2>0$ for $0<p<\sqrt{3}/2$, the positive term $\ph_2^2$ in
the preceding formula cannot vanish, proving the 1st part of the
proposition.

We recall from \eqref{klein:integr2} that
\begin{equation}\label{E1:expr}
E_1=\varphi_1^2+(\varphi_1\varphi_2'-
\varphi_2\varphi_1')^2/3-(\varphi_1\varphi_0'-\varphi_0\varphi_1')^2
\end{equation} Since $E_1<0$ for $0<p<\sqrt{3}/2$, the negative
term $-(\varphi_1\varphi_0'-\varphi_0\varphi_1')^2$ in the
preceding formula cannot vanish, proving the 2nd part of the
proposition.

Finally, if $\varphi_2=\pm 1$ then $\varphi_0=\varphi_1=0$,
contradicting the fact that $\ph_1\ph_0'-\ph_0\ph_1'\neq 0$. \qed

\begin{corollary}\label{cor:1}
For $0<p<\sqrt{3}/2$, the solutions of the system \eqref{syst01}
``rotate'' around the origin in the $(\ph_0,\ph_1)$-plane.  In
other words, if we introduce polar coordinates in \eqref{syst01},
the angle will be monotone increasing.
\end{corollary}

{\bf Proof.} The angle in polar coordinates in the
$(\ph_0,\ph_1)$-plane is given by $\theta=\arctan(\ph_1/\ph_0)$,
and
$$
\theta'=(\varphi_0\varphi_1'-\varphi_1\varphi_0')/(\ph_0^2+\ph_1^2).
$$
Proposition \ref{prop:1} now implies that $\theta'\neq 0$.\qed

Using Corollary \ref{cor:1} we conclude that the condition A in
section \ref{sec:min} (i.e. that $\ph_0,\ph_1$ both have two zeros in the
period) means that a periodic orbit should make {\em one}
turn around the origin.  The periodic solution corresponding to
$p=\sqrt{3/8}$ does exactly that (the orbit in that case is the
ellipse $10\ph_1^2+8\ph_0^2=5$).

\subsection{Intersection angle}\label{sec:symmetry}
Consider the first (for $y>0$) intersection of the trajectory on
the $(\ph_0,\ph_1)$-plane with the $\ph_0$-axis. Let $y(p)$ be the
intersection point, $\alpha(p)$ be the angle of the intersection.
In this section we establish
\begin{prop}\label{prop:symm1}
If $p$ corresponds to an extremal metric, $\alpha(p)=\pi/2$, or,
equivalently, $\ph_0'(y(p))=0$.
\end{prop}

{\bf Proof.} We know from \eqref{func:1}
$$
\ph_0(-y)=\ph_0(y),\ph_1(-y)=-\ph_1(y),
$$
i.e. that the solution for $y>0$ and the solution for $y<0$ are
symmetric with respect to the $\ph_0$-axis.

Assume now that for some $p>0$ the system \eqref{syst01} has a
periodic solution with period $a(p)$.  We have
$$
\ph_0(y(p))=\ph_0(-y(p)),\ \ \ph_1(y(p))=-\ph_1(-y(p))=0.
$$
The periodicity condition together with the condition A imply that
$y(p)=a(p)/2$ and that
$$
\ph_0(y(p)+t)=\ph_0(-y(p)+t), \ \ \ph_1(y(p)+t)=\ph_1(-y(p)+t).
$$
But since $\ph_0$ is an even function, we also have
$$
\ph_0(y(p)+t)=\ph_0(-y(p)+t)=\ph_0(y(p)-t).
$$
The last equality implies $\ph_0'(y(p))=0$.  \qed

\subsection{Ruling out the interval $(0,\sqrt{3}/2)$ numerically}
To rule out the interval $(0,\sqrt{3}/2)$ we use the following
\begin{conj}
\label{conj2} The angle $\alpha(p)$ is a monotone function for $p\in
(0, \sqrt{3}/2)$.
\end{conj}

We check this conjecture numerically for $p\in (\delta,
\sqrt{3}/2-\delta)$ for small $\delta>0$ using {\it Mathematica}.

We have included the values of $\cot(\alpha(p))$ for $999$ values
of $p$, $p=\frac{\sqrt{3}j}{2000},1\leq j\leq 999$.  Those values were
computed using a {\em Mathematica} program. The differential
equation solved by that program is obtained by first writing the
system of two second order differential equations for the
variables $(\psi,\theta)$ corresponding to the change of variables
$$
\ph_2=\cos\psi,\ph_1=\sin\psi\sin\theta,\ph_0=\sin\psi\cos\theta,
$$
then rewriting that system using $\theta$ as an independent
variable (we can do that for $0<p<\sqrt{3}/2$ by Corollary
\ref{cor:1}).

\noindent{The} results are shown below:

\vskip 5pt \hbox to \hsize { \hfill \epsfysize=2in
\epsffile{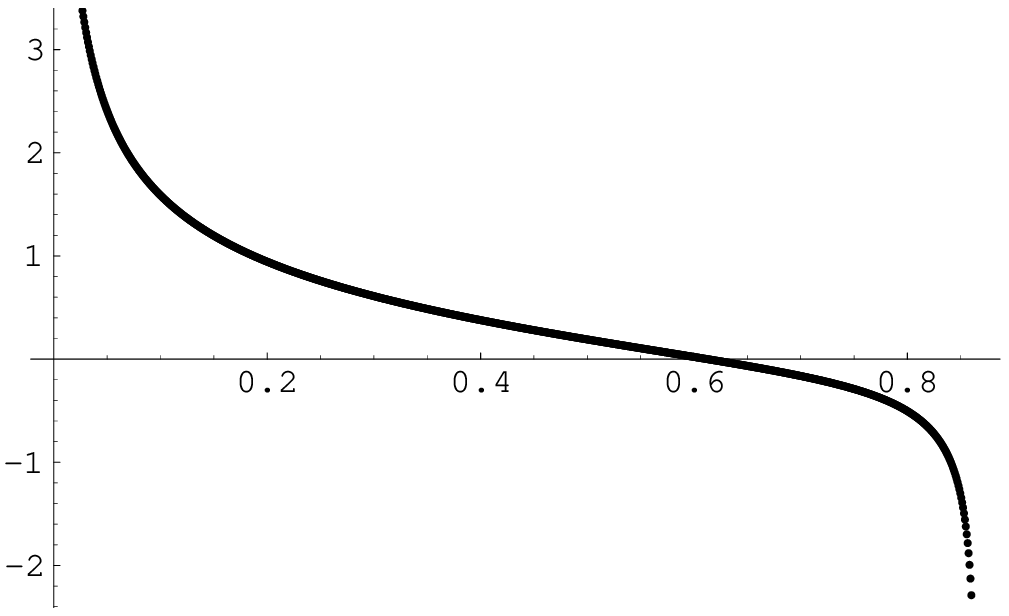} \hfill } \hbox to \hsize { \hfill \hbox to
0pt
{\hss {\em Graph of $\cot(\alpha(p))$ for $0<p<\sqrt{3}/2$.\\
The only zero occurs for $p=\sqrt{\frac{3}{8}}\approx
0.612$.}\hss} \hfill } \vskip 9pt Clearly, since $\cot(\alpha(p))$
is monotone, the same is true for $\alpha(p)$.

We next prove the following
\begin{prop}
Conjecture \ref{conj2} implies Conjecture \ref{conj}.
\end{prop}
\noindent {\bf Proof.} Since the angle $\alpha(p)$ is monotone, it
takes the value $\alpha(p)=\pi/2$ only once. This happens exactly
for $p=\sqrt{3}/8,$ so by Proposition \ref{prop:symm1}, the only
extremal metric on the interval $(0,\sqrt{3}/2)$ is the metric
$g_0$.  \qed

\subsection{Difficulties in proving Conjecture \ref{conj}}
One would like to have a computer-assisted proof of Conjecture
\ref{conj2}, or of a weaker statement (that still implies
Conjecture \ref{conj}) that the conclusion of Proposition
\ref{prop:symm1} only holds for $p=\sqrt{3/8}$.  The main obstacle
to finding a rigorous (even a computer-assisted) proof seems to be
that the systems (\ref{syst01}) and (\ref{syst12}) are lacking
stability, and therefore estimates for the dependence of the
solutions on the initial conditions are very rough.  Accordingly,
one has to make numerical measurements with the step not $10^{-3}$
as in the previous section, but a dozen of orders of magnitude
smaller. Otherwise it seems impossible to control the behavior of
the solutions between the two measurements. However, such
precision seems to be beyond the reach of existing numerical
software.

A similar difficulty occurs near the ends of the interval $(0,
\sqrt{3}/2)$. It can be shown  that
$$
\lim_{p\to 0}y(p)=\lim_{p\to\sqrt{3}/2} y(p)=\infty,
$$
in fact that $y(p)\to\infty$ as $c|\ln p|$ for an explicit
constant $c$. Consequently, if the system \eqref{syst01} for $p>0$
has a periodic solution with period $a(p)$, then
$$
\lim_{p\to 0}a(p)=\lim_{p\to\sqrt{3}/2} a(p)=\infty.
$$
One can also show that there exists an explicit $M>0$ such that for
any $a>M$ and for any metric $g_a$ on $\KK$ with the conformal class
$a$, we have (in the notation of section \ref{sec:main}),
$$
\lambda_1(g_a) {\rm Area}(\KK^2, g_a) < 12\pi E(2\sqrt{2}/3).
$$

Altogether this implies the existence of  a computable constant
$\delta>0$ such that an extremal metric for $\lambda_1 \area$
cannot be attained for $p<\delta$ and $p>\sqrt{3}/2-\delta$.
However, the value of $\delta$ we could obtain is way too small
for being useful in a computer-assisted proof.

{\bf Acknowledgements.}  The authors would like to thank Vestislav
Apostolov,  Maxim Braverman, Alex Ivrii, Niky Kamran, Franz Pedit,
Leonid Pol\-te\-ro\-vich, John Toth, Dimiter Vassilev and
Thomas Zamojski for useful discussions.

A part of this research was conducted
while the third author was visiting the Max Planck Institute for Mathematics
in Bonn. He would like to thank the Institute for its hospitality
and an excellent research atmosphere.

\small

{\bf Note added to the proofs.} Right before the publication of
 this paper we learned that Conjecture \ref{conj} was settled in
\cite{EGJ}. We also thank Hugues Lapointe  for correcting several
inaccuracies in the original version of the present paper.

\end{document}